\documentclass[11pt]{article}
\usepackage{amsmath,amssymb, epsfig}

\newtheorem{theorem}{Theorem}

\newcommand{\R}{{\mathbb R}}

\newcommand{\E}{{\mathbb{E}}}
\newcommand{\Q}{{\mathbb{Q}}}

\catcode`@=11
\@addtoreset{equation}{section}
\catcode`@=12
\begin{document}
\title{Parking on a Random Tree}

\author{H.G.Dehling\footnote{Ruhr-Universit\"at Bochum,
Fakult\"at f\"ur Mathematik,
Universit\"atsstra\ss e 150,
44780 Bochum, Germany,
\texttt{herold.dehling@rub.de}, 
\texttt{http://www.ruhr-uni-bochum.de/ffm/Lehrstuehle/Lehrstuhl-XII/dehling.html}
}, S.R.Fleurke\footnote{Agentschap Telecom, Postbus 450,  
9700 AL Groningen, 
The Netherlands, \texttt{sjoert.fleurke@at-ez.nl}
},\, and 
C. K\"ulske
\footnote{
University of Groningen, 
Department of Mathematics and Computing Sciences, 
Blauwborgje 3,   
9747 AC Groningen, 
The Netherlands,
\texttt{kuelske@math.rug.nl}, 
\texttt{ http://www.math.rug.nl/$\sim$kuelske/ }}
}

\maketitle
\begin{abstract}
Consider an infinite tree with random degrees, i.i.d. over the
sites,  with a prescribed probability distribution with generating
function $G(s)$. We consider the following variation of R\'{e}nyi's
parking problem, alternatively called blocking RSA:
at every vertex of the tree a particle (or ``car'')
arrives with rate one. The particle sticks to the vertex whenever
the vertex and all of its nearest neighbors are not occupied yet.

We provide an explicit expression for the so-called parking
constant in terms of the generating function. That is, the occupation probability, averaged over
dynamics and the probability distribution of the random trees converges in the
large-time limit to $(1-\alpha^2)/2$ with $\int_{\alpha}^1\frac{x dx
}{G(x)}=1$.
\end{abstract}


\vspace{10mm} {\it Key--Words:}
Car parking problem, Random tree, Random sequential adsorption, Particle systems.  

\section{Introduction}

In the classical car parking problem considered by R\'{e}nyi
\cite{Renyi}, one dimensional cars with unit length appear one by
one with their midpoints uniformly distributed over an interval.
A car is parked unless it intersects with one or more previously
placed cars. The process stops, when there is no further possibility
of placing a car. In the discrete version of R\'{e}nyi's car
parking problem, cars of length $2$ try to park at their midpoints
randomly on the integers. Now a car can be parked if the distance
of its midpoint to all other midpoints of already parked cars is two or more.
It is well known that the probability that a given site is occupied
by the midpoint of a car converges to $(1-e^{-2})/2$ in the fully
parked state, that is when time tends to infinity \cite{Co62,
Hemmer}.
 In higher dimensions  a model
of this type will in general not be solvable and more complicated
behavior is expected to occur.
In physical language, this process is called (blocking) RSA (random sequential adsorption),
the motivation being that particles are deposited onto a surface \cite{Co62, Penrose, Evans, Gouet}.
For a nice review article over various  {\it packing} problems (of which the {\it
parking } problem can also be seen as a special example), see
\cite{Cha07}.

In this note we consider the problem of parking on a {\it random}
tree, where particles at the vertices appear with an exponential
waiting time, and parking is prohibited when there is a car-midpoint
already at the site or at a neighboring site {\it on the tree}. We
stress that the tree itself is unchanged under the dynamics. We
solve the model and calculate the occupation probability as a
function of time for random trees with i.i.d. degrees, averaged over
the distribution of trees. In particular we recover the regular tree
result  \cite{Penrose}. We believe that the random tree model has
appeal for two reasons: The first is that  it provides a step to the
analysis of the process on random networks, motivated e.g. by the
study of communication networks. Indeed, many of the random graph
distributions proposed for its description which have found recent
interest  allow for a local random tree approximation \cite{Do03,
Du07, Ho07}. The second reason why it is nice is just the
mathematical simplicity of the result which should not go unnoted.

The proof is based on
an analysis of the ordinary differential equation of
the occupation probability at a fixed site of a fixed realization of the
random tree. This is like  \cite{Penrose, Hemmer}, but we give a self-contained
exposition.
Conditioning on a non-arrival of a car at this site
the r.h.s. of this differential equation
allows for a factorization in terms of occupation probabilities on rooted trees.
Again a differential equation for these quantities can be derived.
While it cannot be expected to be solvable for  a particular realization,
the averaging over the tree distribution still allows for closed-form
expressions, as we shall see.

\subsection{Random Trees with Independent Identically Distributed Degrees}

The precise definition of the model is as follows. We consider a
random tree with vertices $i$ and degree at the site $i$ given by
$D_i$. We choose $D_i$ to be independent random variables with the
same distribution $ \mathbb{Q}$ given by
\begin{equation}\begin{split}
 \mathbb{Q}(D_i=k)=a_k
\end{split}\end{equation}
on the integers starting from $2$. The latter requirement ensures
that we have no open ends with probability one. We denote the
generating function of the distribution by
\begin{equation}\begin{split}
 G(s)=\sum_{k=2}^\infty a_k s^k
\end{split}\end{equation}

We will denote the expected value with respect to this probability
distribution by the same symbol $\Q$.

\subsection{The Dynamics}

Fix a tree with vertex set $V$. For any such tree
we will define a Markov jump process on the occupations numbers
$n=(n_i)_{i\in V}\in \Omega=\{0,1\}^V$.  Here $n_i$ denotes the occupation number
of vertex $i$, meaning that
\[ n_{i} = \left\{ \begin{array}{rl}
0 & \mbox{if vertex $i$ is vacant } \\
1 & \mbox{if vertex $i$ is occupied }  \end{array} \right. \]

The dynamics of the process is defined in terms of the
generator which is given by the r.h.s. of the differential equation
\begin{equation}\begin{split}
\frac{d}{dt} \mathbb{E}^{n} f(n(t)) = \sum_{k} \left[ f(n^k) - f(n)
\right] r_k(n) \label{Generator}
\end{split} \end{equation}
with
\[ n^{k}_i = \left\{ \begin{array}{rl}
n_i & \mbox{if $k \neq i$ } \\
1 & \mbox{if $k = i$ }  \end{array} \right. \] and
\begin{equation}\begin{split}
r_k(n)&=\prod_{i:d(i,k)\leq 1}n_{i}^c, \cr
n^c _i &= 1 - n_i,
\end{split}\end{equation}
and where $d$ denotes the distance on the tree.

This generator defines a Markov jump process on the infinite graph
by standard theory \cite{Liggett}, such that  (\ref{Generator})
holds for any local function $f: \Omega\rightarrow\R$. Here $\E^{n}$
denotes the expected value with respect to the process, started at
the initial configuration $n$.

This generator describes the parking of cars at all possible sites $k$
at rate $1$ given a configuration which is given
by $n$.  This parking is possible whenever the site and its nearest neighbors are
vacant, that is $r_k(n)=1$.
Since at most one car can be parked at a given site, it suffices to consider
the time of the first arrival of a car at this site, for all sites.

\subsection{Results}

We provide an explicit integral formula for the density of occupied
sites at time $t$, averaged with respect to dynamics $\E$ and tree
distribution $\Q$.

\begin{theorem}
\begin{equation}\begin{split}
\Q\E(n_0(t))= \frac{1 - \alpha^2 \bigl(1-e^{-t}\bigr)}{2}
\end{split}\end{equation}
where the function $\alpha(u)$ is defined by the equation
\begin{equation}\begin{split}\label{1.6}
\int_{\alpha(u)}^1\frac{x dx }{G(x)}=u
\end{split}\end{equation}
In particular, the occupation probability, averaged over
dynamics and the probability distribution of the random trees converges in the
large-time limit to $(1-\alpha^2)/2$ with $\int_{\alpha}^1\frac{x dx
}{G(x)}=1$.
\end{theorem}
\bigskip

{\bf Remark: } Let us specialize to the deterministic case of a
regular tree with $D\geq 2$ nearest neighbors. We obtain
\[
\E(n_0(t)) = \left\{
  \begin{array}{ll}
\frac{1}{2}\bigl(1-e^{-2(1-e^{-t})}\bigr)& \mbox{if  } D=2  \\[2mm]
\frac{1}{2}\bigl(1 -  (1+(D-2)(1-e^{-t}))^{-\frac{2}{D-2}}\bigr)
 & \mbox{if } D>2
 \end{array}
 \right.
 \]
 and recover in this way the known results on the integers and general regular trees
\cite{Co62, Hemmer, Penrose}.

\section{Proof}
Consider a fixed realization of the random tree and look at the time-evolution
given by (\ref{Generator}).
We use a short notation for the probability that a given set of sites $A$ is empty at time $t$,
\begin{equation}\begin{split}
C_t(A)&:=\E(\prod_{i\in A}n^c_{i}(t))\cr
\end{split}\end{equation}
Similarly we write for the conditional probability that a given set of sites $A$ is empty, conditional on the event
that on another set of sites $B$ no particles have arrived yet:
\begin{equation}\begin{split}
C_t(A|B)&:=\E(\prod_{i\in A}n^c_{i}(n)| T_{j}>t \text{ for } j\in B )\cr
\end{split}\end{equation}
Of course, if no particle has arrived at a site it is empty, but the
converse is not true. Here $T_j$ denotes the arrival time of the
particle at site $j$. The $T_j$'s are i.i.d. exponentially
distributed with expected value equal to one. We also write in short
$C_t(i_1,\dots,i_k)\equiv C_t(\{i_1,\dots,i_k\})$ for
$A=\{i_1,\dots,i_k\}$ etc.

As the underlying tree is random, also
these correlation functions are random variables with respect to the
distribution $\Q$. For any  realization we can write the
dynamics (\ref{Generator}) as
\begin{equation}\begin{split}\label{2.2}
-\frac{d}{dt}C_t(A)&=C_t(\overline{A})\cr
\end{split}\end{equation}
where $\overline A=\{ i\in G| d(i,A)\leq 1  \}$. Similarly we have for
the dynamics of the conditional correlation function the differential equation
\begin{equation}\begin{split}\label{2.3}
-\frac{d}{dt}C_t(A| B)&=C_t(\overline{A}| B)\cr
&=C_t(\overline{A}\cap B^c| B)   \cr
\end{split}\end{equation}
The first equation follows from restricting the generator to the set
$B^c$ (since in $B$ no particles have arrived yet). The second
equation expresses the fact that when no particles have arrived in
the set $B$, it implies that $\bar A \cap B$ is unoccupied.

Fix an arbitrary vertex. Let us call this vertex $0$.
Starting from the dynamics for the correlation functions (\ref{2.2}) we have
\begin{equation}\begin{split}
-\frac{d}{dt}C_t(0)&=C_t(\{i,d(i,0)\leq 1\})\cr
&=\E\Bigl(\prod_{i,d(i,0)\leq 1}n^c_{i}(t)
\,\,1_{T_0 >t } \Bigr)
+\E\Bigl(\prod_{i,d(i,0)\leq 1}n^c_{i}(t)
\,\, 1_{T_0 \leq t } \Bigr)\cr
\end{split}\end{equation}
We notice that the second term vanishes since an arrival of a particle at $0$
would imply that the sites $\{i,d(i,0)\leq 1\}$ cannot all be empty.
Writing the first term as a conditional probability we get
\begin{equation}\begin{split}\label{2.6}
-\frac{d}{dt}C_t(0)&=C_t(\{i,d(i,0)\leq 1\} | 0)\, e^{-t}\cr
\end{split}\end{equation}
Conditional on the event that no
particle has arrived at the site $0$, the dynamics for the
occupation numbers on the branches of the tree that are emerging from
$0$ is independent. Consequently the correlations factorize into a
product and we have
\begin{equation}\begin{split}\label{2.7}
C_t(\{i,d(i,0)\leq 1\}| 0)=\prod_{i,d(i,0)=1}C_t(i | 0)\cr
\end{split}\end{equation}
where the $C_t(i | 0)$ is the non-occupation probability of a site $i$
that is adjacent to the root on a rooted tree, assuming that no
particle has arrived at the root.

\begin{figure}
\setlength{\unitlength}{1mm}
\begin{center}
\begin{picture}(50,60)(0,-7)
\multiput(5,20)(0,5){5}{\circle*{2}} %
\put(5,20){\line(1,1){10}} \put(5,25){\line(2,1){10}}
\put(5,30){\line(1,0){10}} \put(5,35){\line(2,-1){10}}
\put(5,40){\line(1,-1){10}} \put(15,30){\circle*{2}}
\put(24,30){\line(-1,0){10}} \put(24,30){\line(-1,1){10}}
\put(26,30){\line(1,1){10}}  \put(26,30){\line(1,0){10}}
\put(25,29){\line(0,-1){8}} \put(25,30){\circle{2}}
\put(25,30){\circle{1}}
\put(25,20){\circle{2}} \put(25,20){\circle{1}}%
\put(24,20){\line(-1,-1){10}} \put(25,19){\line(0,-1){10}}
\put(26,20){\line(1,-1){10}} \put(14,10){\circle*{2}}
\put(25,10){\circle*{2}} \put(36,10){\circle*{2}}
\put(13,8){$\underbrace{~~~~~~~~~~~~~~~~~~~}$} \put(20,2){$D_1-1$}
\put(35,30){\circle*{2}} \put(35,30){\line(1,1){10}}
\put(35,30){\line(1,-1){10}} \put(45,40){\circle*{2}}
\put(45,20){\circle*{2}}
\put(36,40){\circle*{2}} \put(36,40){\line(0,1){10}}
\put(36,50){\circle*{2}} %
\put(14,40){\circle*{2}} \put(14,40){\line(0,1){10}}
\put(14,40){\line(1,1){10}} \put(14,50){\circle*{2}}
\put(24,50){\circle*{2}}
\put(25,32){0} \put(26,21){1}
\end{picture}
\caption{Fragment of a tree where sites 0 and 1 are empty and the
occupancy of the other sites is unknown.} \label{fig:tree}
\end{center}
\end{figure}
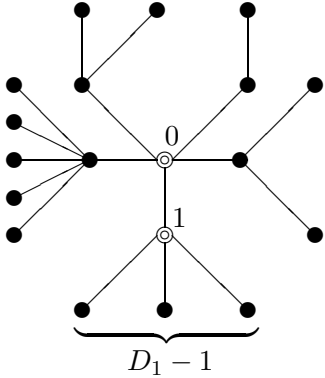

For a particular realization of the underlying tree, these functions
will in general differ. Let us however now take the expected value
over the tree distribution $\Q$. We can decompose this average into
an average over $D_0$ (the number of neighbors of the central vertex
$0$) and a conditional average over the remaining branches, rooted
at $0$. Conditioning the number of branches emerging from zero to be
equal to $k$, the correlation functions are independent random
variables with respect to $\Q$ and so we obtain
\begin{equation}\begin{split}
\Q \Bigr( C_t(\{i,d(i,0)=1\}| 0)\Bigl | D_0=k\Bigl) =\Bigl( \Q
\Bigr( C_t(1 |  0) \Bigl)\Bigl) ^k
\end{split}\end{equation}
where we denoted by the site $1$ one of the nearest neighbors of
$0$, see Fig. \ref{fig:tree}.

From here, (\ref{2.6}) and (\ref{2.7}) we have
\begin{equation}\begin{split}\label{2.9}
-\frac{d}{dt}\Q (C_t(0))
&=\sum_{k=2}^\infty a_k\Bigl( \Q C_t(1 | 0)\Bigr) ^k e^{-t}\cr
&=G\bigl( \Q C_t(1 | 0)\bigr) e^{-t}
\end{split}\end{equation}
Integrating this equation gives us the desired function $\Q
(C_t(0))$ once the ``half-infinite function''  $\Q C_t(1 | 0)$ is
known. So, it remains to solve a problem on a rooted tree,
conditioned that no particle has arrived at the root, averaged over
$\Q$.

Indeed, it turns out to be possible to derive a closed-form
differential equation for this object, as we show now. Note at first
that it suffices to consider the event that no particle has arrived
at $1$ and write
\begin{equation}\begin{split} \label{2.10}
C_t(\{i,d(i,1)\leq 1\}| 0)=C_t(\{i,d(i,1)\leq 1\}| 0,1)e^{-t}\cr
\end{split}\end{equation}
Since the root $0$ and the site $1$ itself are empty under the conditioning we have
\begin{equation}\begin{split}
C_t(\{i,d(i,1)\leq 1\}| 0,1)&= C_t(\{i\neq 0,d(i,1)=1\}| 0,1)\cr
&= \prod_{i\neq 0,d(i,1)=1}C_t(i | 0,1)
\end{split}\end{equation}
Employing the dynamics  (\ref{2.3}), using (\ref{2.10}) and
averaging over $\Q$ we thus have
\begin{equation}\begin{split}\label{2.12}
-\frac{d}{dt}\Q C_t(1 | 0)&=\Q C_t(\{i,d(i,1)\leq 1\}| 0)\cr
&=\Q \Bigl( \prod_{i\neq 0,d(i,1)=1}C_t(i | 0,1) \Bigr) e^{-t}\cr
\end{split}\end{equation}
Using conditional independence over the branches with respect to the
tree distribution $\Q$ and noting that $\Q C_t(i | 0,1) = \Q C_t(1
| 0)$ for the $i$'s appearing in the above product to get
\begin{equation}\begin{split}\label{2.13}
\Q \Bigl( \prod_{i\neq 0,d(i,1)=1}C_t(i | 0,1) \Bigr)
=\sum_{k=2}^\infty a_k (\Q C_t(1 | 0))^{k-1}
\end{split}\end{equation}
we finally obtain with (\ref{2.12}) and (\ref{2.13})
the closed-form differential equation
\begin{equation}\begin{split}
-\frac{d}{dt}\Q C_t(1 | 0)&=\frac{G\bigl( \Q C_t(1 | 0)\bigr)}{\Q C_t(1 | 0)} e^{-t}\cr
\end{split}\end{equation}
Its solution $\Q C_{t}(1| 0)$ with initial condition $\Q C_{t=0}(1| 0)=1$ is given by the integral
\begin{equation}\begin{split}
\int_{\Q C_t(1| 0)}^{1}\frac{y dy }{G(y)}=1-e^{-t}
\end{split}\end{equation}

Let us substitute this solution into the equation for the filling probability (\ref{2.9}).
This then implies the explicit formula  for the $\Q$-average of the filling density of a random tree with generating function $G$
of the form
\begin{equation}\begin{split}
 \Q \E (n_0(t))&=\int_{0}^t  G\bigl( y(s)\bigr) e^{-s} ds
\end{split}\end{equation}
where $y(s)$ is defined by the equation
\begin{equation}\begin{split} \label{2.17}
\int_{y(s)}^1\frac{x dx }{G(x)}=1-e^{-s}
\end{split}\end{equation}
Finally, with the substitution $u=1-e^{-s}$ and recalling
(\ref{1.6}) we see that
\begin{equation}\begin{split}
\int_{0}^t G(y(s))e^{-s} &= \int_{0}^{1-e^{-t}} G(y(u))du =- \int_1 ^{\alpha(1-e^{-t})} y \, dy \cr
&= \frac{1}{2}(1-\alpha^2(1-e^{-t}))
\end{split}\end{equation}
where the second step follows from (\ref{2.17}). This finishes the proof.

\end{document}